\documentclass[preprint, 5p]{elsarticle}

\usepackage{datetime}
\newdateformat{UKvardate}{%
\THEDAY\ \monthname[\THEMONTH] \THEYEAR}
\UKvardate

\usepackage{hyperref}
\usepackage{subcaption}
\usepackage{xcolor}
\usepackage{todonotes}

\usepackage{booktabs}
\usepackage{float}

\usepackage{tabularx}
\newcolumntype{C}{>{\centering\arraybackslash}X}
\usepackage{multirow}

\usepackage{mathtools}
\usepackage{amsmath}
\usepackage{amssymb}
\usepackage{amsthm}
\usepackage{marvosym} 

\newtheorem{theorem}{Theorem}
\newtheorem{lemma}[theorem]{Lemma}

\theoremstyle{definition}

\theoremstyle{remark}

\usepackage{etoolbox}
\makeatletter
\patchcmd{\@thm}{\thm@headpunct{.}}{\thm@headpunct{}}{}{}
\makeatother

\makeatletter
\@ifclassloaded{beamer}{%
    \usepackage[algoruled, noline, noend]{algorithm2e}
}{%
    \usepackage[algoruled, noline, linesnumbered, noend]{algorithm2e}
}
\makeatother

\usepackage{interval}
\intervalconfig{soft open fences}

\usepackage{bbm}
\newcommand{\bb}[1]{\mathbb{#1}}
\newcommand{\cl}[1]{\mathcal{#1}}

\newcommand{\diag}[1]{\mathop{}\!\mathrm{diag}\left(#1\right)}
\newcommand{\lspan}[1]{\mathop{}\!\mathrm{span}\left(#1\right)}
\newcommand{\trace}[1]{\mathop{}\!\mathrm{tr}\left(#1\right)}

\newcommand{\iu}{{i\mkern1mu}}

\newcommand{\set}[1]{\left\{#1\right\}}

\DeclareMathOperator*{\argmax}{arg\,max}

\DeclarePairedDelimiterXPP{\abs}[1]{}{\lvert}{\rvert}{}{{#1}}
\DeclarePairedDelimiterXPP{\norm}[2]{}{\lVert}{\rVert}{_{#2}}{#1}
\DeclarePairedDelimiterXPP{\inner}[3]{}{\langle}{\rangle}{_{#3}}{#1, #2}

\usepackage{enumitem}

\title{Randomized Methods for Kernelized DMD}
\author[1]{Peter Oehme\corref{cor}}
\affiliation[1]{%
    organization={Institute of Mathematics, EPFL},%
    country={Switzerland},%
}
\cortext[cor]{Corresponding author: \emph{oehme.pb@gmail.com}}

\begin{document}
    \begin{abstract}
        Dynamic Mode Decomposition (DMD) is a data-driven method related to Koopman operator theory that extracts information about dominant dynamics from data snapshots. In this paper we examine techniques to accelerate the application of DMD to large-scale data sets with an eye on randomized techniques. Randomized techniques exploit low-rank matrix approximations at a much smaller computational cost, therefore permitting the use of increased data set sizes. In particular, we propose the application of the RPCholesky algorithm in the setting of kernelized DMD (KDMD). This algorithm relies on adaptive randomized sampling to approximate positive semidefinite kernel matrices and provides better stability guarantees than previously implemented randomized methods for KDMD. Differences between existing competitive randomized techniques and our proposed implementation are discussed with a focus on numerical stability and tradeoff between exploration and exploitation of information obtained from data. The efficacy of this new combination of algorithms is demonstrated on well-established benchmark problems from DMD literature increasing in problem dimension.
    \end{abstract}
    \maketitle

    \section{Introduction}

    Dynamic mode decomposition (DMD) is a numerical method that aims to extract the dominant dynamics from sequential snapshots of some underlying nonlinear dynamical system in the form of DMD modes~\cite{Schmid2010}. The DMD modes are derived from the left eigenvectors of the finite-dimensional projections of the so-called Koopman operator --- the operator mapping nonlinear observables of the system's state forward in time. For this purpose, multiple algorithms have been developped~\cite{Tu2014, Colbrook2024} many of which focus on particular attributes such as structure preservation or the induction of sparsity. One of the constraints of dynamic mode decomposition is that datasets with large state dimensions $N$ and large numbers of snapshots $M$ require a lot of computational effort. This is a clear bottleneck of DMD and multiple methods have been proposed aiming to remedy this constraint.

    A prominent technique is kernelized DMD (KDMD)~\cite{Williams2015}, where one applies the kernel trick~\cite{Schoelkopf2000, Baudat2001} to the input data, thus reducing the effects of the large state dimension $N$. The resulting algorithm, however, requires the solution of two $M \times M$ eigenproblems and therefore does not provide a good scaling with the number of snapshots. On the other hand, Randomized Fourier features (RFF)~\cite{Rahimi2007} provide an approach to only approximate the kernel function involved in KDMD, and then use these features in the extended DMD (EDMD) algorithm~\cite{DeGennaro2019} as the nonlinear observables. RFF is based on an analytic relation between the kernel function and a probability measure, allowing for a randomized approach. The principal issue of RFF is the slow convergence of the randomized approximation~\cite{Sutherland2015}, but it has been observed that RFF still reesults in good nonlinear feature maps.
    
    Alternatively, one may use algebraic methods for kernel approximations, in particular for symmetric positive semidefinite kernels via the oASIS algorithm~\cite{Patel2015, Litzinger2018}. This algorithm computes approximate eigenpairs of the kernel matrix by selecting columns from it by diagonal pivoting~\cite{Higham2002} and incrementally updating a kernel matrix approximant via the Schur complement. This often works quite well, however there exist counterexamples where diagonal pivoting fails to return an accurate approximation~\cite{Chen2024}.

    Moreover, the computation of the inverse matrix in the oASIS algorithm is mathematically equivalent to the computation of a pivoted partial Cholesky factorization~\cite[p.~24]{Zhang2005}, however the direct implementation of oASIS is numerically unstable in contrast to the numerical pivoted partial Cholesky factorization.

    The main purpose of this paper is to combine the adaptive random sampling of the RPCholesky algorithm~\cite{Chen2024} with KDMD to address the stability problems and the drawbacks of the diagonal pivoting of the oASIS algorithm. We provide observations suggesting that the RPCholesky algorithm can lead to more stable approximations of the DMD modes when ordering the resulting DMD modes by means of data-driven residuals in the sense of~\cite{Drmac2018}. We demonstrate the efficacy of the randomized KDMD algorithm and compare it to relevant competing algorithms for multiple numerical experiments.

    The paper is organized as follows: In Section~\ref{sec:dmd} we give a short overview of the DMD algorithms relevant to this paper, focusing on how to derive data-driven residuals from the algorithm's outputs. In Section~\ref{sec:randomized-kdmd} we express different options of how randomized methods can be used with KDMD, emphasizing the usage of the RPCholesky algorithm, and analyze the computational complexities of the mentioned algorithms with respect to the dimensions of the input data. Lastly, we present numerical experiments in Section~\ref{sec:num-exp}.

    \section{Fundamentals of DMD}\label{sec:dmd}

    We now explain the background of DMD relevant to this paper.

    \subsection{SVD-Based DMD}

    In~\cite{Schmid2010} the classical SVD-based DMD algorithm was introduced. Given the data matrices $X$ and $Y \in \bb{C}^{N \times M}$
    \begin{equation*}
        X = [x^{(1)}, x^{(2)}, \dots, x^{(M)}], \quad Y = [x^{(2)}, x^{(3)}, \dots, x^{(M + 1)}],
    \end{equation*}
    we compute the SVD $X = U \Sigma V^*$, truncate $U, \Sigma$ and $V$ to the numerical rank $R > 0$, and define
    \begin{equation}\label{eq:svd-dmd-K}
        K = U^* Y V \Sigma^\dagger \in \bb{C}^{R \times R}.
    \end{equation}
    Then, the DMD modes $\xi_i = U w_i \in \bb{C}^N$ computed from the right eigenpairs $(\lambda_i, w_i)$ of $K$ span a low-dimensional subspace that captures dynamics from the data trajectory $x^{(1)}, x^{(2)}, \dots, x^{(M)}, x^{(M + 1)}$. To identify the dominant modes in this subspace we use the data-driven residual $r(i)$ in the sense of~\cite{Drmac2018} given by
    \begin{equation}\label{eq:dmd-residual}
        r(i) = \norm{Y V \Sigma^\dagger w_i - \lambda_i \xi_i}{2}
    \end{equation}
    and select the modes $\xi_i$ with indices $i$ minimizing the residual $r(i)$.

    \subsection{EDMD}

    The EDMD algorithm developped in~\cite{Williams2015a, Kawahara2016} is similar to the SVD-based DMD algorithm, however EDMD includes the use of (nonlinear) observables $\Psi \colon \bb{C}^N \to \bb{C}^{1 \times L}$. Here, the data matrices become $\Psi_X$ and $\Psi_Y \in \bb{C}^{M \times L}$
    \begin{equation*}
        \Psi_X = \begin{bmatrix}
            \Psi(x^{(1)}) \\
            \Psi(x^{(2)}) \\
            \vdots \\
            \Psi(x^{(M)})
        \end{bmatrix},\quad \Psi_Y = \begin{bmatrix}
            \Psi(x^{(2)}) \\
            \Psi(x^{(3)}) \\
            \vdots \\
            \Psi(x^{(M + 1)})
        \end{bmatrix}.
    \end{equation*}
    For the computation of the DMD modes we now define the $L \times L$ matrix $K$ using the SVD $\Psi_X = U \Sigma V^*$ as
    \begin{equation}\label{eq:edmd-K}
        K = \Psi_X^\dagger \Psi_Y,
    \end{equation}
    and from the left eigenpairs $(\lambda_i, w_i)$ and the right eigenpairs $(\lambda_i, \eta_i)$ of $K$ we compute the DMD modes $\xi_i$ and the eigenfunctions $\varphi_i(z) \in \bb{C}$, following the derivations in~\cite{Williams2015a, Williams2015}, by scaling $w_i^* \eta_j = \delta_{ij}$ we have
    \begin{equation*}
        \xi_i = X V \Sigma^\dagger U^* w_i,\quad \varphi_i(z) = \Psi(z) \eta_i,
    \end{equation*}
    respectively. From $\varphi_i(X) = \Psi_X \eta_i \in \bb{C}^M$ we obtain the reconstruction
    \begin{equation}\label{eq:reconstruction}
        \hat{X} = \sum_i \xi_i \varphi_i(X)^\top.
    \end{equation}
    Because of the similarities between~\eqref{eq:svd-dmd-K} and~\eqref{eq:edmd-K}, we use the analogous data-driven residual in the sense of~\cite{Drmac2018}
    \begin{equation*}
        r(i) = \norm{\Psi_Y^* \xi_i - \lambda_i w_i}{2}
    \end{equation*}
    to select the dominant modes and eigenfunctions.

    \subsection{KDMD}

    Lastly, KDMD~\cite{Williams2015} implicitly chooses the observables $\Psi$ by selecting a kernel function $k(x, y) = \Psi(y) \Psi(x)^*$, motivated by the kernel trick in machine learning applications~\cite{Schoelkopf2000, Baudat2001}. With the Gramians $G$ and $A \in \bb{C}^{M \times M}$
    \begin{equation*}
        G = k(X, X) = \Psi_X \Psi_X^*,\quad A = k(X, Y) = \Psi_Y \Psi_X^*
    \end{equation*}
    we compute the eigendecomposition $G = U \Sigma^2 U^*$, and construct the $M \times M$ matrix
    \begin{equation}\label{eq:kdmd-K}
        \hat{K} = \Sigma^\dagger U^* A U \Sigma^\dagger.
    \end{equation}
    Again, from the left eigenpairs $(\lambda_i, w_i)$ and the right eigenpairs $(\lambda_i, \eta_i)$ of $K$ we compute the DMD modes $\xi_i$ and the eigenfunctions $\varphi_i(z)$ by evaluating
    \begin{align*}
        \xi_i &= X U \Sigma^\dagger w_i, \\
        \varphi_i(z) &= \Psi(z) V \eta_i = \Psi(z) \Psi_X^* U \Sigma^\dagger \eta_i = k(X, z) U \Sigma^\dagger \eta_i,
    \end{align*}
    respectively. It particularly holds that $\varphi_i(X) = U \Sigma \eta_i$. These eigenpairs of $K$ from~\eqref{eq:edmd-K} and $\hat{K}$ from~\eqref{eq:kdmd-K} are related by the following statement.

    \begin{lemma}[{\cite[Prop.~1]{Williams2015}}]
        Let $U \Sigma V^*$ be the SVD of $\Psi_X$. Then, $\lambda \neq 0$ and $\eta$ are a right eigenpair of $\hat{K}$ from~\eqref{eq:kdmd-K} if and only if $\lambda$ and $\hat{\eta} = V \eta$ are a right eigenpair of $K$ from~\eqref{eq:edmd-K}. In an analogous way it follows that $\lambda \neq 0$ and $\eta$ are a left eigenpair of $\hat{K}$ from~\eqref{eq:kdmd-K} if and only if $\lambda$ and $\hat{\eta} = V \eta$ are a left eigenpair of $K$ from~\eqref{eq:edmd-K}.
    \end{lemma}

    We require a different method to compute a data-driven residual in the sense of~\cite{Drmac2018}, because we can no longer apply~\eqref{eq:dmd-residual}. We concentrate on the quality of $w_i$ as a left eigenvector and see that
    \begin{align*}
        K^* \hat{w}_i &= \lambda_i \hat{w}_i = \lambda_i V w_i = \lambda_i \Psi_X^* U \Sigma^\dagger w_i, \text{ and} \\
        K^* \hat{w}_i &= \Psi_Y^* {(\Psi_X^*)}^\dagger \hat{w}_i = \Psi_Y^* U \Sigma^\dagger V^* \hat{w}_i = \Psi_Y^* U \Sigma^\dagger w_i,
    \end{align*}
    whereafter we take the difference of both terms and multiply by $\Psi_X$ from the left to arrive at
    \begin{equation}\label{eq:kdmd-residual}
        r(i) = \norm{A^* U \Sigma^\dagger w_i - \overline{\lambda}_i G^* U \Sigma^\dagger w_i}{2}
    \end{equation}
    as a suitable candidate for a data-driven residual in the KDMD algorithm.

    Once more we select the dominant DMD modes and eigenfunctions by choosing the indices $i$ with the smallest residuals $r(i)$. In total, we write the KDMD algorithm in Algorithm~\ref{alg:kdmd} in pseudocode.

    \begin{algorithm}[ht!]
        \caption{KDMD}\label{alg:kdmd}
        \KwData{Data matrices $X, Y \in \bb{C}^{N \times M}$, kernel $k$}
        \KwOut{Eigenvalues $\Lambda = [\lambda_1, \lambda_2, \dots, \lambda_R]$, evaluated eigenfunctions $\Phi(X) \in \bb{C}^{M \times R}$, DMD modes $\Xi \in \bb{C}^{R \times N}$, and data-driven residuals $r(i)$}
        Compute the Gramians $G = k(X, X), A = k(X, Y) \in \bb{C}^{M \times M}$\;
        Compute the eigendecomposition $G = U \Sigma^2 U^*$\;\label{algline:kdmd-eprob1}
        Truncate to the numerical rank $R$ $U = U[:, : R] \in \bb{C}^{M \times R}, \Sigma = \Sigma[: R, : R] \in \bb{R}^{R \times R}$\;
        Form $\hat{K} = \Sigma^\dagger U^* A U \Sigma^\dagger \in \bb{C}^{R \times R}$\;
        Compute the eigenpairs $W, V, \Lambda = \texttt{eig}(\hat{K})$ with left eigenvectors $W$ and right eigenvectors $V$\;\label{algline:kdmd-eprob2}
        Scale $W$ such that $w_i^* v_j = \delta_{ij}$, the Kronecker delta\;
        Compute the eigenfunctions $\Phi(X) = U \Sigma V \in \bb{C}^{M \times R}$\;
        Compute the DMD modes $\Xi = X U \Sigma^\dagger W \in \bb{C}^{R \times N}$\;
        Compute the residuals $r(i)$ from~\eqref{eq:kdmd-residual}\;
    \end{algorithm}

    The choice of a kernel for the KDMD algorithm is a difficult problem and such a choice depends on the underlying data~\cite{Williams2015}. To keep it simple, we restrain ourselves to the use of radial basis functions, in particular the Gaussian kernel.

    \section{Randomized sampling for KDMD}\label{sec:randomized-kdmd}

    We observe that the dominating cost of the KDMD algorithm are the two $M \times M$ eigenproblems (lines~\ref{algline:kdmd-eprob1} and~\ref{algline:kdmd-eprob2}). These bottlenecks prohibit the application of KDMD to datasets with long snapshot trajectories $M \gg 0$. There exist some proposed solutions to this issue, of which we will particularly mention random Fourier features and the Nyström approximation.

    \subsection{Random Fourier Features (RFF)}

    The seminal paper~\cite{Rahimi2007} makes use of the fact that for symmetric, translation-invariant and positive semidefinite kernels $k$ Bochner's theorem states that $k$ can be expressed as the Fourier transform of a probability measure. Sampling $z_1, z_2, \dots, z_S, S > 0$, from this measure~\cite{DeGennaro2019} proposes to compute the observables $\psi_\ell(x) = \exp{(\iu \langle x, z_\ell \rangle)}$ for $\ell = 1, 2, \dots, S$ which approximate the kernel $k$ by the sample average
    \begin{equation*}
        k(x, y) \approx \frac{1}{S} \sum\limits_{i = 1}^S \psi_i(x) \overline{\psi_i(y)}.
    \end{equation*}
    Afterwards, one uses the observables $\psi_i(x)$ in the EDMD algorithm. This means that RFF does not explicitly aim at providing a speedup of the KDMD algorithm, but rather tries to find a good choice of observables for the EDMD algorithm.

    Unfortunately, the convergence of this approximation is relatively slow with the expected maximum error decaying with a rate of $\cl{O}(1/S)$, cf.~\cite{Sutherland2015}. Additionally, the kernel's underlying probability measure is often unavailable, hence requiring either empirical estimations of the measure or the use of specific kernels such as the Gaussian kernel for which it is known that the underlying measure is again Gaussian, but with inverse covariance.

    \subsection{Diagonal Pivoting and oASIS}\label{subsec:diagonal-pivoting}

    The oASIS algorithm~\cite{Patel2015} computes the Nyström approximation~\cite{Williams2000}
    \begin{equation*}
        G \approx \tilde{G} = G[I, :] {G[I, I]}^{-1} G[:, I]^*,
    \end{equation*}
    where $I = \set{i_1, i_2, \dots, i_S}$ is a set of $S$ suitably chosen pivot elements. These pivot elements are chosen one after another such that
    \begin{align*}
        i_{k + 1} = \argmax\limits_{j = 1, 2, \dots, M} \diag{G - C_k {G[I_k, I_k]}^{-1} C_k^*}[j],
    \end{align*}
    where $C_k = G[:, I_k]$ is the column-selected matrix and the set $I_k = \set{i_1, i_2, \dots, i_k}, k < S$, is an intermediate pivot selection. The authors then successively update the matrices ${G[I_k, I_k]}^{-1}$ and $C_k {G[I_k, I_k]}^{-1}$ by use of the Schur complement to finally extract the eigenpairs of $\tilde{G}$ as approximations of the eigenpairs of the original kernel matrix $G$.

    The pivoting rule in~\cite{Patel2015} turns out to be equivalent to diagonal pivoting for the pivoted Cholesky algorithm~\cite[Section~10.3]{Higham2002} because the matrix $G - C_k {G[I_k, I_k]}^{-1} C_k^*$ coincides with the residual $A - L_k L_k^*$ of the pivoted Cholesky algorithm and because
    \begin{align*}
        i_{k + 1} &= \argmax\limits_{j = 1, 2, \dots, M}\ (d - \mathrm{colsum}(C_k \odot R_k))[j] \\
                  &= \argmax\limits_{j = 1, 2, \dots, M}\ G[j, j] - \sum\limits_{\ell} C_k[j, \ell] R_k[j, \ell] \\
                  &= \argmax\limits_{j = 1, 2, \dots, M}\ G[j, j] - \sum\limits_{\ell} C_k[j, \ell] ({G[I_k, I_k]}^{-1} C_k)[\ell, j] \\
                  &= \argmax\limits_{j = 1, 2, \dots, M}\ \diag{G - C_k {G[I_k, I_k]}^{-1} C_k^*}[j]
    \end{align*}
    for $d = \diag{G}$ and $R_k = G[:, I_k] {G[I_k, I_k]}^{-*}$. Unlike the Schur complement of the oASIS algorithm, however, the pivoted Cholesky algorithm is numerically stable. The motivation for the diagonal pivoting is to avoid numerical breakdown whilst computing the square root of the main diagonal element $G[i_{k + 1}, i_{k + 1}]$ for positive semidefinite matrices. While this works well if one executes the Cholesky algorithm until termination, there exist counterexamples where the greedy pivot selection and early termination lead to suboptimal partial approximations~\cite{Chen2024}. 

    We thus replace the oASIS algorithm which has already been used to speed up the identification of molecular dynamics in an EDMD-like framework, cf.~\cite{Litzinger2018} by the RPCholesky algorithm~\cite{Chen2024}, hence benefitting from both the increased numerical stability and the more robust randomized pivot selection. At the same time we preserve the advantage that replacing $G$ by $\tilde{G}$ reduces the $M \times M$ eigenproblem in line~\ref{algline:kdmd-eprob1} of Algorithm~\ref{alg:kdmd} to an $S \times S$ eigenproblem.

    \subsection{Randomly Pivoted Partial Cholesky Factorization}

    Instead of the diagonal pivoting choice from Section~\ref{subsec:diagonal-pivoting}, we utilize the randomized selection proposed in~\cite{Chen2024}. This means that we again choose pivots $I = \set{i_1, i_2, \dots, i_S}, S > 0$, but instead of greedily choosing an index from the largest elements of the remainder's diagonal elements we balance this greedy exploitation with randomized exploration of the smaller diagonal elements. In particular, the next pivot $i_{k + 1}$ is chosen according to the probability distribution given by the remainder matrix $G_k = G - C_k {G[I_k, I_k]}^{-1} C_k^*$'s diagonal elements relative weights w.r.t.\ to each other as
    \begin{equation*}
        \bb{P}(i_{k + 1} = j) = \frac{\diag{G_k}[j]}{\trace{G_k}}.
    \end{equation*}

    It has been shown that this procedure known as \emph{RPCholesky}~\cite{Chen2024} produces a good approximation of the optimal rank-$S$ approximation of $G$ with high probability. Unlike the diagonal pivoting the random selection strikes a balance between exploitation of maximum incoherence and random exploration of unknown data. It has also been demonstrated that this behaviour can outperform greedy selections depending on the kernel and the data set. We use the resulting factorization $G \approx F F^*$ to replace the eigenvector computation in line~\ref{algline:kdmd-eprob1} of Algorithm~\ref{alg:kdmd} by the eigendecomposition
    \begin{equation*}
        F^* F = \hat{U} \Sigma^2 \hat{U}^* 
    \end{equation*}
    and the orthonormalization $U = \mathrm{orth}(F \hat{U})$. This change has the additional side effect that the matrix $\hat{K}$ will also only be of shape $R \times R$ with the numerical rank $R < S$. We summarize our proposed alteration in Algorithm~\ref{alg:rand-kdmd}.

    \begin{algorithm}[ht!]
        \caption{KDMD with randomized sampling}\label{alg:rand-kdmd}
        \KwData{Data matrices $X, Y \in \bb{C}^{N \times M}$, kernel $k$}
        \KwOut{Eigenvalues $\Lambda = [\lambda_1, \lambda_2, \dots, \lambda_R]$, evaluated eigenfunctions $\Phi(X) \in \bb{C}^{M \times R}$, DMD modes $\Xi \in \bb{C}^{R \times N}$, and data-driven residuals $r(i)$}
        Compute the Gramian $A = k(X, Y)$ and run RPCholesky $F, I = \texttt{rpcholesky}(X, k)$\;
        Compute the eigendecomposition $F^\top F = U \Sigma^2 U^*$ and orthonormalize $U = \mathrm{orth}(F \hat{U})$\;
        Truncate to the numerican rank $R$ $U = U[:, : R] \in \bb{C}^{M \times R}, \Sigma = \Sigma[: R, : R] \in \bb{R}^{R \times R}$\;
        Form $\hat{K} = \Sigma^\dagger U^* A U \Sigma^\dagger \in \bb{C}^{R \times R}$\;
        Compute the eigenpairs $W, V, \Lambda = \texttt{eig}(\hat{K})$ with left eigenvectors $W$ and right eigenvectors $V$\;
        Scale $W$ such that $w_i^* v_j = \delta_{ij}$, the Kronecker delta\;
        Compute the eigenfunctions $\Phi(X) = U \Sigma V \in \bb{C}^{M \times R}$\;
        Compute the DMD modes $\Xi = X U \Sigma^\dagger W \in \bb{C}^{R \times N}$\;
        Compute the residuals $r(i)$ from~\eqref{eq:rkdmd-residual}\;
    \end{algorithm}

    Besides Algorithm~\ref{alg:rand-kdmd} we also propose a data-driven residual in the sense of~\cite{Drmac2018} for the resulting DMD modes. Here, we reuse the information used in the RPCholesky execution, namely $\hat{G}[:, I]$ and the selection of pivots $I$. We define the residual $r(i)$ of the $i$-th eigenvalue DMD mode pair $(\lambda_i, \xi_i) = (\lambda_i, U \Sigma^\dagger W)$ as
    \begin{gather}\label{eq:rkdmd-residual}
        \begin{aligned}
            r(i) &= \norm{A[:, I]^* \xi_i - \overline{\lambda} G[:, I]^* \xi_i}{2} \\
                 &= \norm{A[:, I]^* U \Sigma^\dagger w_i - \overline{\lambda} G[:, I]^* U \Sigma^\dagger w_i}{2}.
        \end{aligned}
    \end{gather}
    We use this residual to sort the DMD modes from Algorithm~\ref{alg:rand-kdmd} and select the dominant modes with the lowest residuals.

    In addition to Algorithm~\ref{alg:rand-kdmd} we highlight that for sequential data $X$ and $Y$ like we use here we can further simplify the precalculation of the matrix $A$. Because $X$ and $Y$ are subamtrices of the matrix
    \begin{equation*}
        D = [x^{(1)}, x^{(2)}, \dots, x^{(M)}, x^{(M + 1)}]
    \end{equation*}
    we can assemble $k(D, D)$ and extract $G$ and $A$ from the upper left block and the upper right block, respectively. Thus, we do not need to assemble $A$ separately, but can directly use $\hat{F}, \hat{I} = \texttt{rpcholesky}(k(D, D))$, and $G \approx F F^\top$ and $A \approx F B^\top$ for $F = \hat{F}[:-1, :]$ and $B = \hat{F}[:1, :]$. This way we save more computational effort in the preparatory stages of kernel matrix assembly.

    \subsection{Computational Complexity of Randomized KDMD}\label{subsec:runtime}

    We quickly list the runtimes of the algorithms we covered so far in Figure~\ref{tab:runtimes}. We split the analysis of the complexities into three parts: Firstly, the assembly of the matrix of observables or the kernel matrices, secondly, the assembly of the matrices $K$ and $\hat{K}$ depending on the algorithm, and lastly, the assembly of the final DMD modes. For the EDMD algorithm the runtime depends on the number of observables $L$, however a bad set of observables can lead to $L \in \cl{O}(N)$.

    Notably, for the KDMD algorithm using the pivoted partial Cholesky factorization we can compute the decomposition of the kernel matrix for the entire data matrix instead of only $X$ or $Y$ directly because we have assumed sequential data. This means that the computation of the decomposition only requires $S$ evaluations of colums of the entire kernel matrix, resulting in a complexity of $\cl{O}(S M N)$, and an additional $\cl{O}(S^2 M)$ operations to run the algorithm itself. This way we can avoid the expensive complete assembly of $G$ and $A$ entirely, saving on computational costs. The KDMD algorithm using the pivoted partial Cholesky algorithm thus uses $\cl{O}(S M N + S^2 M + S^3)$ operations in total.
    \begin{table*}[ht!]
        \centering
        \begin{tabular}{c c c c}
            Algorithm & Observables/Kernel Matrix & Assembly of $K$/$\hat{K}$ & DMD Modes \\
            \hline
            EDMD~\cite{DeGennaro2019, Colbrook2023} & $\cl{O}(L M N)$ & $\cl{O}(L^2 M + L^3)$ & $\cl{O}(L M N + L^2 N + L^3)$ \\
            KDMD~\cite{DeGennaro2019} & $\cl{O}(M^2 N)$ & $\cl{O}(M^3)$ & $\cl{O}(M^2 N + M^3)$ \\
            RFF EDMD~\cite{DeGennaro2019} & $\cl{O}(S M N)$ & $\cl{O}(S^2 M + S^3)$ & $\cl{O}(S M N + S^2 N + S^3)$ \\
            Algorithm~\ref{alg:rand-kdmd}~\cite{Chen2024} & $\cl{O}(S M N + S^2 M)$ & $\cl{O}(S^2 M + S^3)$ & $\cl{O}(S M N + S^3)$
        \end{tabular}
        \caption{Computational complexity of the algorithms covered. The oASIS algorithm because it has the same complexity as Algorithm~\ref{alg:rand-kdmd}.}\label{tab:runtimes}
    \end{table*}

    \section{Numerical Experiments}\label{sec:num-exp}

    The code for these experiments can be found on GitHub under \href{https://github.com/peoe/randomized-kdmd}{github.com/peoe/randomized-kdmd}. We run the experiments on a MacBook Air 2022 with an M2 processor and 16 GB of RAM.

    In all experiments, whenever we refer to KDMD or one of its randomized variants, we will use the Gaussian kernel with the scaling parameter $\sigma = 1 / N$. This choice reflects the default choice for the implementation of the radial basis function kernels in \texttt{scikit-learn}.

    \begin{figure*}[ht!]
        \centering
        \subfloat[DMD]{\includegraphics[width=.45\textwidth]{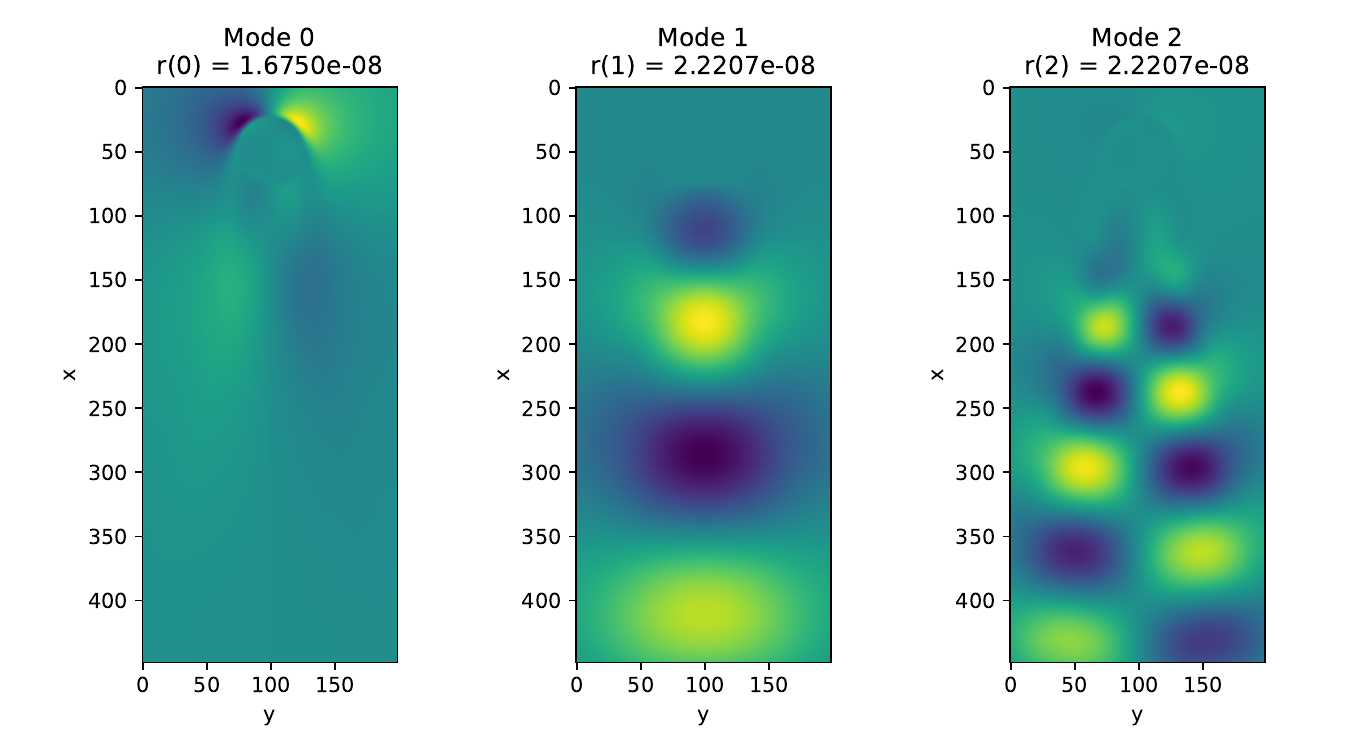}}
        \subfloat[KDMD]{\includegraphics[width=.45\textwidth]{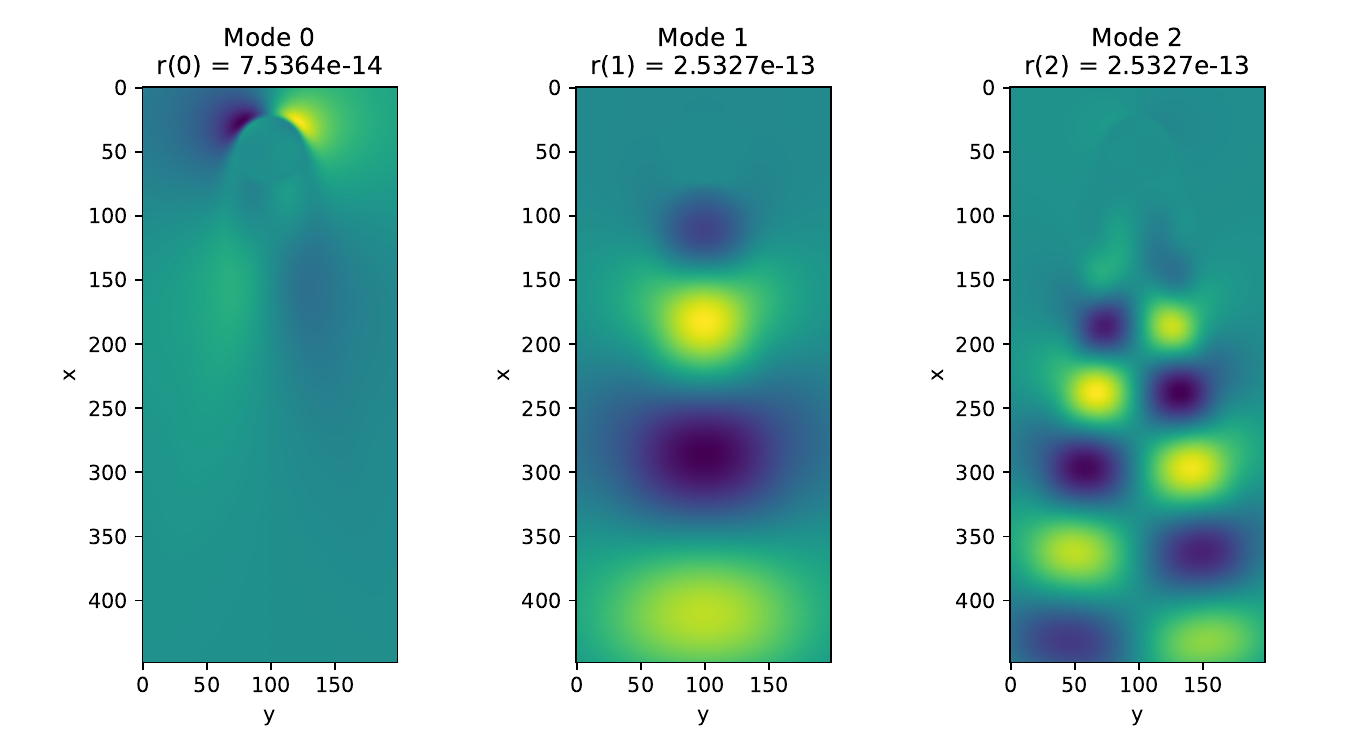}} \\
        \subfloat[RFF]{\includegraphics[width=.45\textwidth]{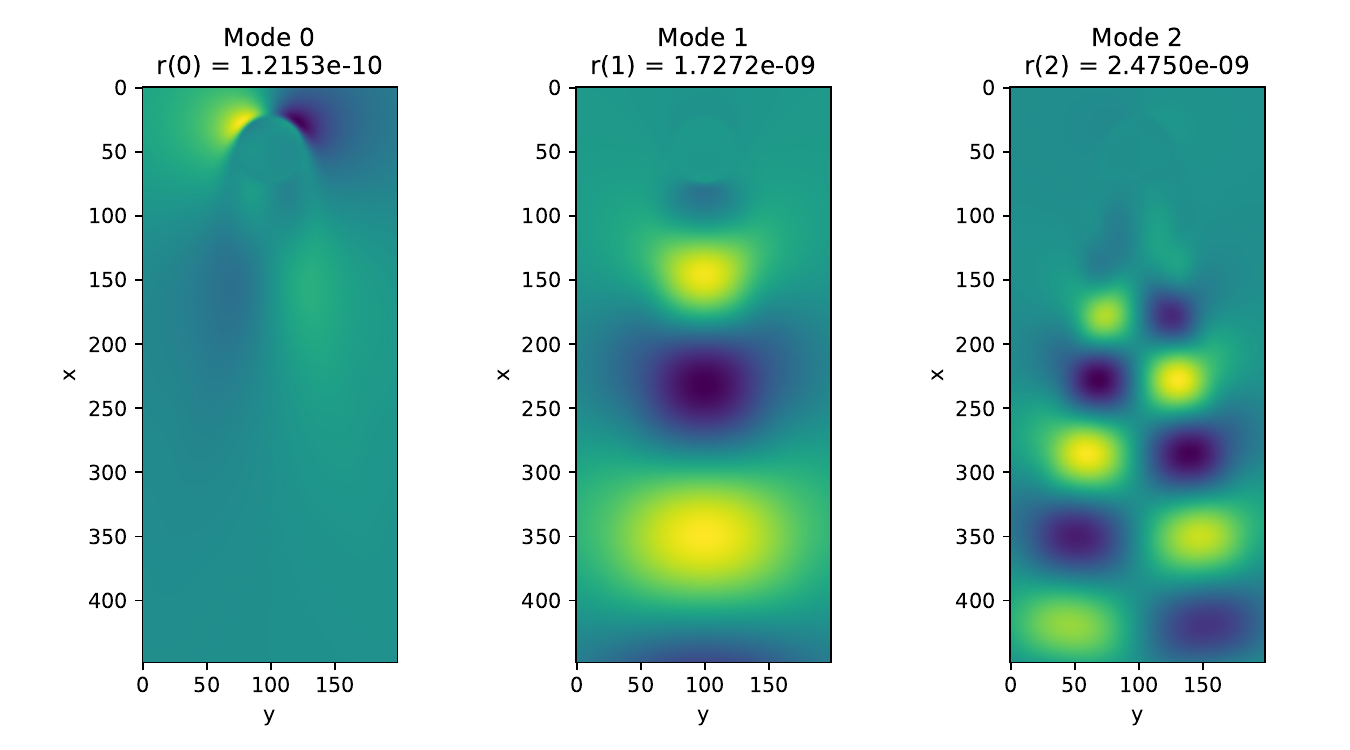}}
        \subfloat[oASIS]{\includegraphics[width=.45\textwidth]{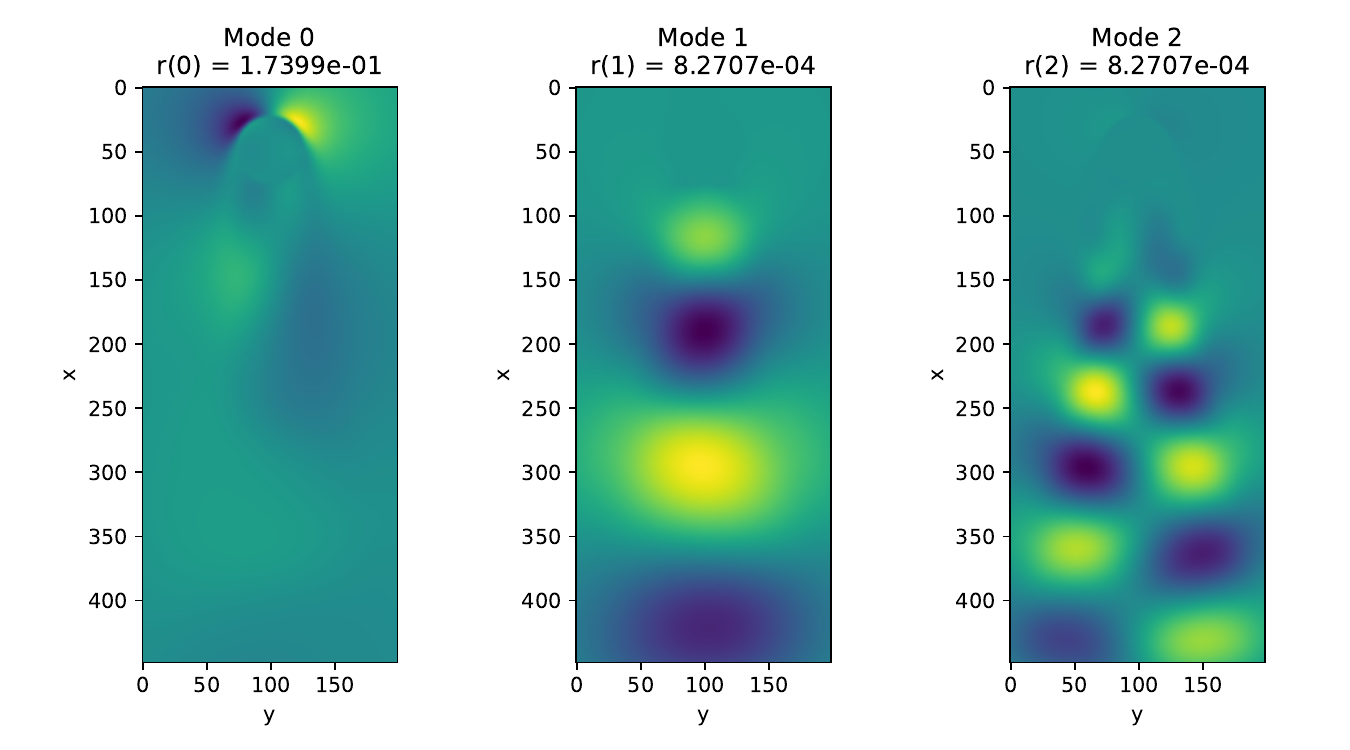}} \\
        \subfloat[Algorithm~\ref{alg:rand-kdmd}]{\includegraphics[width=.45\textwidth]{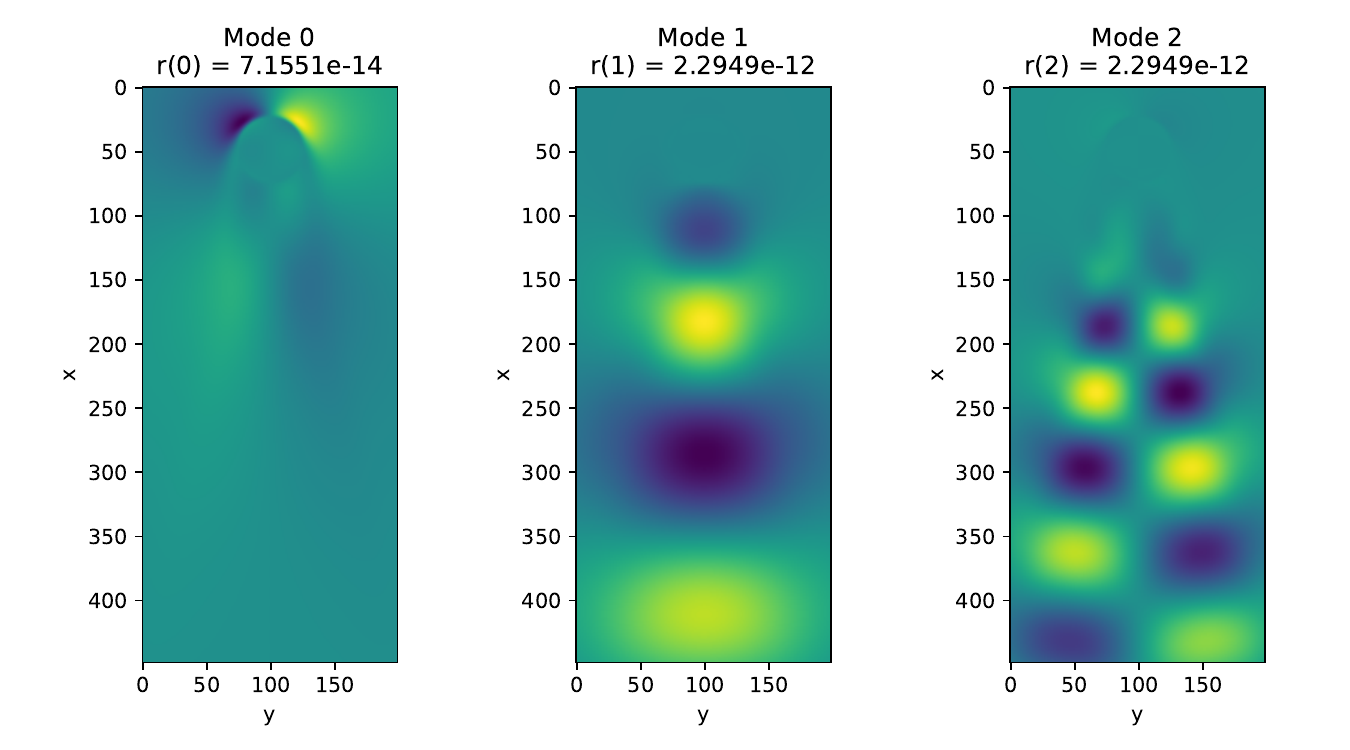}}
        \caption{The real parts of the first three DMD modes for fluid flow around cylinder. The modes were increasingly sorted according to their respective residuals $r(i)$. For the randomized methods RFF and RPCholesky we used $S = 50$ random samples. For RFF, KDMD with RPCholesky, and KDMD with greedily pivoted partial Cholesky the signs are flipped, but the same phenomenon occurs to the corresponding eigenfunctions.}\label{fig:cylinder-flow-modes}
    \end{figure*}

    \subsection{Flow Around a Cylinder}

    We start with a popular experiment in DMD literature, the laminar flow of a fluid around a cylinder at Reynolds number $\mathrm{Re} = 100$. This example is based on the two-dimensional incompressible Navier-Stokes equations
    \begin{align*}
        \frac{\partial v}{\partial t} &= -(v \cdot \nabla) v + \nu \Delta v - \frac{1}{\rho} \nabla p, \\
        \nabla \cdot v &= 0,
    \end{align*}
    where $v$ is the two-dimensional velocity field, $p$ is the pressure of the fluid, $\rho$ is the fluid's density, and $\nu$ is the kinematic viscosity. We only use the pressure data from numerical simulations available in the data provided with the DMD Book~\cite{Kutz2016}. Each of the pressure snapshot trajectories have a shape of $89351 \times 150$, that is $N = 89351$ and $M = 150$. Clearly, this data does not require the RPCholesky decomposition because even standard DMD is sufficient for this example, however this makes it easy to qualitatively compare the modes computed by different DMD algorithms.

    We run SVD-based DMD, default KDMD, RFF EDMD with samples from a Gaussian distribution with $\mu = 0, \sigma=1 / 10$, and the KDMD algorithm with randomized matrix decompositions on the $v$-velocity data. We choose the matrix decompositions from RPCholesky, greedily pivoted partial Cholesky, and Nyström approximation with pivot selection via oASIS. In Figure~\ref{fig:cylinder-flow-modes} we plot the first five modes, omitting two plots which are complex conjugates of the displayed modes.

\begin{figure}[ht!]
        \centering
        \includegraphics[width=.45\textwidth]{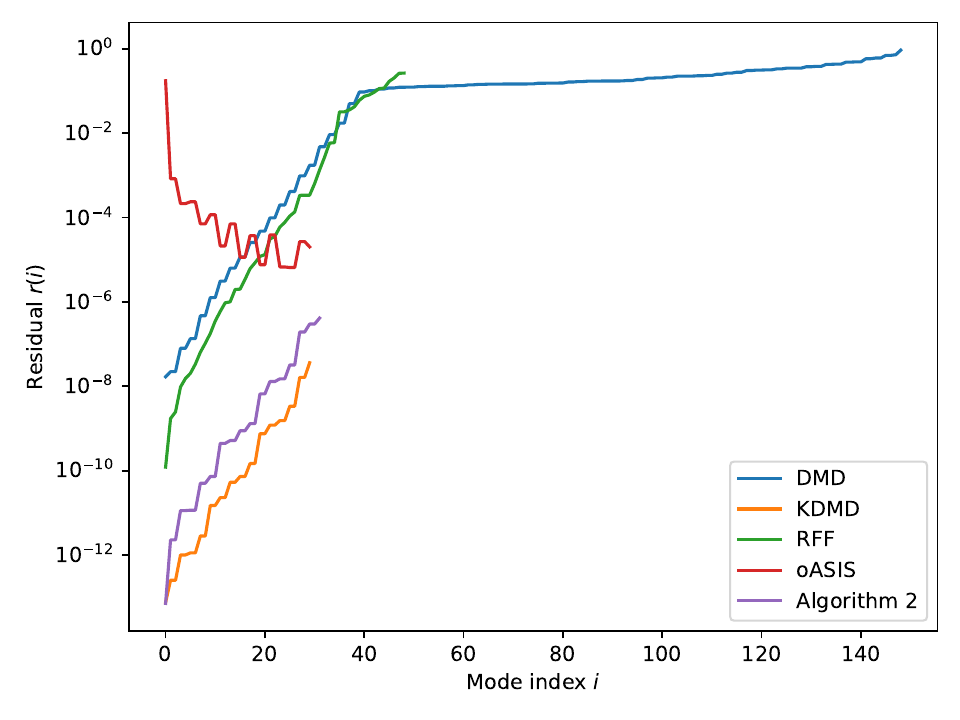}
        \caption{Residuals of the DMD modes for the fluid flow around a cylinder. All randomized methods use $S = 50$ random samples.}\label{fig:cylinder-flow-residuals}
    \end{figure}

    In Figure~\ref{fig:cylinder-flow-residuals} we plot the data-driven residuals in the sense of~\cite{Drmac2018} in comparison to each other. Noticeably, the residual of the Nyström approximation via oASIS pivot selection does not reflect the same behaviour as the residuals of the other methods. Due to the erratic behaviour of the residuals we do not sort the modes increasingly by their residuals for this method, preserving instead the order in which they were computed. Besides this we remark that the modes are quite similar visually.

    \subsection{Unforced Duffing Equation}

    In this experiment we consider a trajectory of the unforced Duffing equation similar to~\cite{Williams2015a}. The governing equation for the Duffing oscillator with two basins of attraction that we will consider is
    \begin{equation}\label{eq:duffing-oscillator}
        \ddot{x} = -\frac{1}{2} \dot{x} - x (-1 + x^2)
    \end{equation}
    We select $700$ random initial points for $(x, \dot{x})$ from a uniform distribution over the square $\interval{-2}{2}^2$ and compute $11$ time steps for each of these initial conditions. We concatenate all these trajectories into a single matrix of size $2 \times 7700$, and repeat the matrix along the first axis to arrive at the final data matrix of size $10000 \times 7700$. This data matrix will probably not contain a single trajectory of the dynamical system~\eqref{eq:duffing-oscillator} that converges in one of the basins of attraction, however it contains enough data that we can identify the overall dynamics and therefore get an idea of the basins of attraction. Additionally, this data matrix is only of rank $2$, however its size lets it serve as a stress test for the KDMD algorithms.

\begin{figure}[ht!]
        \centering
        \includegraphics[width=.45\textwidth]{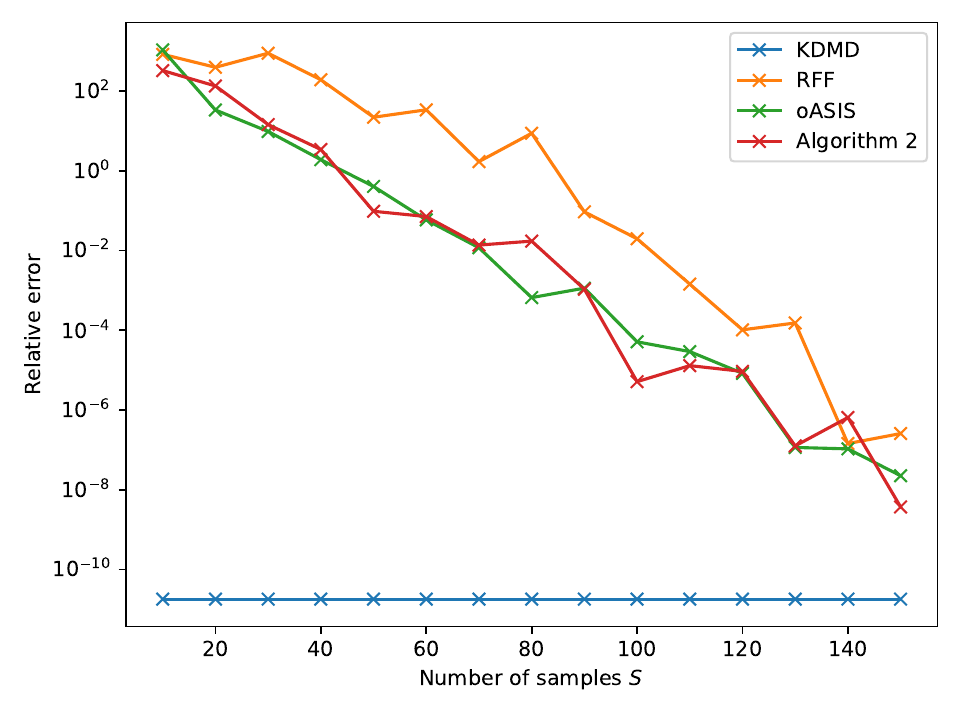}
        \caption{Reconstruction errors for the synthetic Duffing oscillator example. The shape parameters of the Gaussian distribution generating the RFF samples and the Gaussian kernel used in the KDMD algorithms are chosen such that $\sigma = 1/\sqrt{2}$}\label{fig:duffing-reconstruction-errors}
    \end{figure}

    We now compare the relative reconstruction errors
    \begin{equation*}
        \tau = \sum\limits_{m = 1}^M \frac{\norm*{x^{(m)} - \hat{x}^{(m)}}{2}}{\norm{x^{(m)}}{2}}
    \end{equation*}
    for the default KDMD algorithm, KDMD with diagonally pivoted Cholesky factorization, KDMD with RPCholesky, KDMD with oASIS selection for Nyström pivots, and RFF EDMD. For the random Fourier features we choose a Gaussian distribution with the shape parameter $\sigma = 1 / 20$. The reconstruction $\hat{x}$ here is given by~\eqref{eq:reconstruction}
    \begin{equation*}
        \hat{X} = \sum\limits_i \xi_i \varphi_i(X)^*.
    \end{equation*}
    We successively increase the number of samples $S$ and plot the resulting relative errors $\tau$ in Figure~\ref{fig:duffing-reconstruction-errors}.

\begin{figure}[ht!]
        \centering
        \includegraphics[width=.45\textwidth]{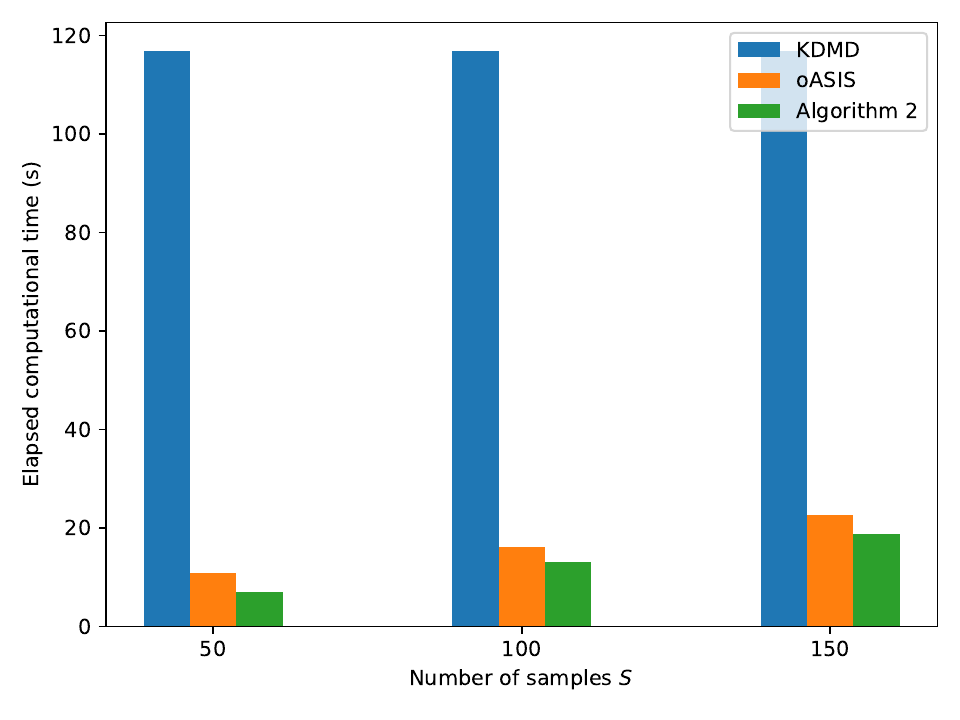}
        \caption{Elapsed computational times for the synthetic Duffing oscillator data set.}\label{fig:duffing-times}
    \end{figure}

    In Figure~\ref{fig:duffing-times} we plot the computational times for increasing ranks for all KDMD variants. We notice that the required computational time for the randomized KDMD algorithms rises with an increasing number of samples $S$. Importantly, all three methods appear to be scaling with the same rate with respect to $S$, confirming the runtime analysis done in Section~\ref{subsec:runtime}.

    \subsection{Sea Surface Temperature Dataset}

    Lastly, we demonstrate the capabilities of the randomized algorithms mentioned in this work in application to a large dataset. To this end, we choose the weekly average NOAA Sea Surface Temperature data set available at \href{https://psl.noaa.gov/data/gridded/data.noaa.oisst.v2.highres.html}{psl.noaa.gov/data/gridded/data.noaa.oisst.v2.highres.html} and previously examined in~\cite{Erichson2019} and~\cite{Askham2018}. In total, it contains $2276$ snapshots of $1036800$ data points, hence resulting in matrices $G$ and $A$ of dimensions $2276 \times 2276$. While this dimension is inferior to the matrix size in the Duffing oscillator example, the evaluation of the kernel matrices $G$ and $A$ without preassembly of the whole matrix is much slower.

    Besides our proposed randomized algorithm~\cite{Erichson2019} promulgated an implementation of SVD-based DMD where they used a randomized SVD algorithm to accelerate the SVD. Unfortunately, this means that one would still need to compute a QR decomposition $Q R = X \Omega, \Omega \in \bb{C}^{M \times S},$ of the large matrix $X \Omega \in \bb{C}^{N \times S}$. While this makes the overall application of the DMD algorithm viable, one still needs to preprocess all the available data by projecting $X$ onto $\lspan{Q}$ by explicitly multiplying $\tilde{X} = Q^* X$. This means that we can easily store $\tilde{X}$ all at once, but the preprocessing nevertheless requires a large amount of time. The advantage of the randomized KDMD algorithms over the randomized SVD-based DMD algorithm is that we only ever require access to a maximum of two snapshot at a time.

\begin{figure}[ht!]
        \centering
        \subfloat[RFF]{\includegraphics[width=.45\textwidth]{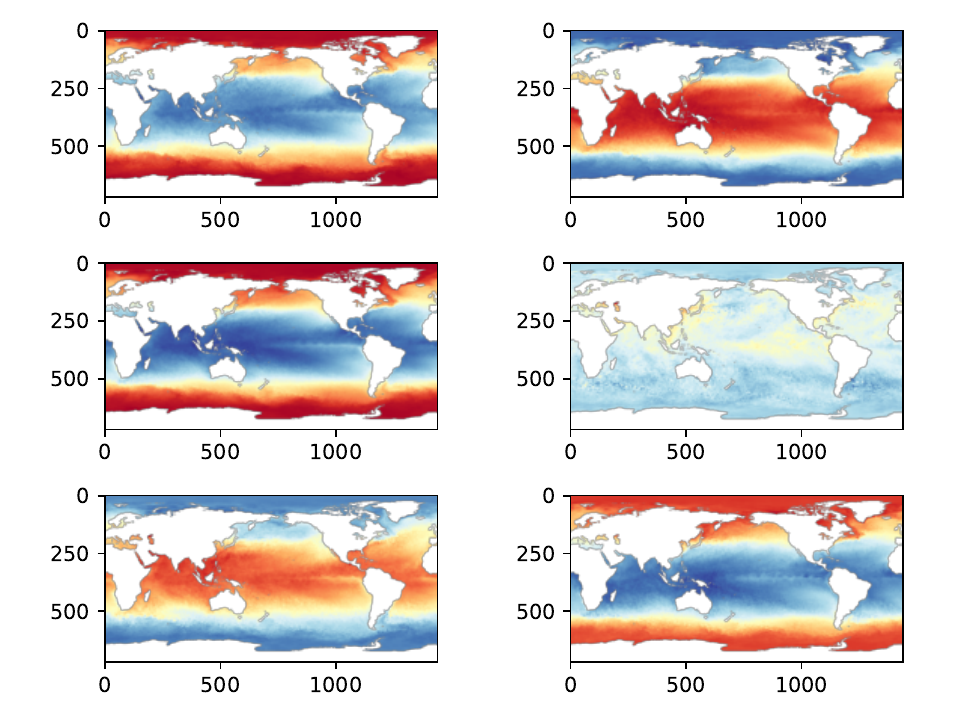}} \\
        \subfloat[Algorithm~\ref{alg:rand-kdmd}]{\includegraphics[width=.45\textwidth]{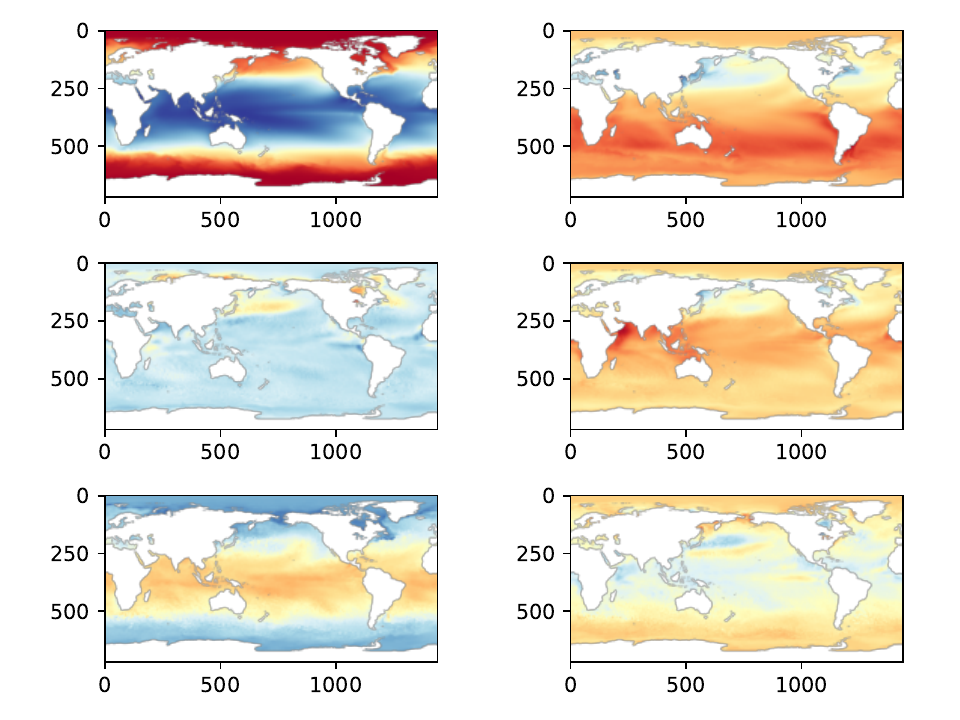}}
        \caption{The first 6 DMD modes for the sea surface temperature data set. Here, both techniques are randomized and use $S = 100$ as the sample size.}\label{fig:sst-modes}
    \end{figure}

    In Figure~\ref{fig:sst-modes} we compare the first 6 modes produced by the KDMD algorithm with RPCholesky with those produced by RFF EDMD. For both computations we chose $S = 100$ samples. The random Fourier features are computed from samples of a Gaussian distribution with shape parameter $\sigma = 1 / 100$. We remark that the modes produced by Algorithm~\ref{alg:rand-kdmd} fulfill the usual assumptions about DMD modes, namely that the first mode appears as a standalone mode whereas all subsequent modes appear in complex conjugate pairs. In contrast, the modes computed by RFF EDMD do not result in this pairwise ordering.

    \section{Conclusion}

    We demonstrated that randomization and partial pivoted Cholesky factorizations are viable tools to accelerate the KDMD algorithm and derived a DMD mode residual formula in the sense of~\cite{Drmac2018} resulting in a good ordering of the computed DMD modes. We showed empirically that the randomization via RPCholesky improves upon the computation based on the oASIS algorithm, in particular because the resulting residuals provide a meaningful ordering of the DMD modes. Additionally, we highlighted that the oASIS pivot selections correspond to those of the diagonally pivoted partial Cholesky factorization. The latter factorization provides the same benefits as the RPCholesky factorization in the fact that it also results in a meaningful residual. Diagonal pivoting, however, may suffer from the drawbacks of only exploiting available information instead of weighted random exploration as RPCholesky chooses.

    \bibliographystyle{elsarticle-num}
    \bibliography{rkdmd}

@InProceedings{Rahimi2007,
  author    = {Rahimi, Ali and Recht, Benjamin},
  booktitle = {Advances in Neural Information Processing Systems},
  title     = {Random Features for Large-Scale Kernel Machines},
  year      = {2007},
  publisher = {Curran Associates, Inc.},
  volume    = {20},
}

@InProceedings{Sutherland2015,
  author    = {Sutherland, Danica J. and Schneider, Jeff},
  booktitle = {Proceedings of the Thirty-First Conference on Uncertainty in Artificial Intelligence},
  title     = {On the error of random {F}ourier features},
  year      = {2015},
  address   = {Arlington, Virginia, USA},
  pages     = {862--871},
  publisher = {AUAI Press},
  series    = {UAI'15},
  isbn      = {9780996643108},
  location  = {Amsterdam, Netherlands},
  numpages  = {10},
}

@Article{DeGennaro2019,
  author     = {DeGennaro, Anthony M. and Urban, Nathan M.},
  journal    = {SIAM J. Sci. Comput.},
  title      = {Scalable extended dynamic mode decomposition using random kernel approximation},
  year       = {2019},
  issn       = {1064-8275,1095-7197},
  number     = {3},
  pages      = {A1482--A1499},
  volume     = {41},
  doi        = {10.1137/17M115414X},
  fjournal   = {SIAM Journal on Scientific Computing},
  mrclass    = {37M10 (62M15 65P99)},
  mrnumber   = {3948243},
  mrreviewer = {Rua\ Murray},
}

@Article{Chen2024,
  author   = {Chen, Yifan and Epperly, Ethan N. and Tropp, Joel A. and Webber, Robert J.},
  journal  = {Comm. Pure Appl. Math.},
  title    = {Randomly pivoted {C}holesky: Practical approximation of a kernel matrix with few entry evaluations},
  year     = {2024},
  doi      = {10.1002/cpa.22234},
  fjournal = {Communications on Pure and Applied Mathematics},
}

@Article{Williams2015,
  author   = {Williams, Matthew O. and Rowley, Clarence W. and Kevrekidis, Ioannis G.},
  journal  = {J. Comput. Dyn.},
  title    = {A kernel-based method for data-driven {K}oopman spectral analysis},
  year     = {2015},
  issn     = {2158-2491,2158-2505},
  number   = {2},
  pages    = {247--265},
  volume   = {2},
  doi      = {10.3934/jcd.2015005},
  mrnumber = {3507895},
}

@Article{Erichson2019,
  author   = {Erichson, N. Benjamin and Mathelin, Lionel and Kutz, J. Nathan and Brunton, Steven L.},
  journal  = {SIAM J. Appl. Dyn. Syst.},
  title    = {Randomized dynamic mode decomposition},
  year     = {2019},
  issn     = {1536-0040},
  number   = {4},
  pages    = {1867--1891},
  volume   = {18},
  doi      = {10.1137/18M1215013},
  fjournal = {SIAM Journal on Applied Dynamical Systems},
  mrnumber = {4021266},
}

@Article{Williams2015a,
  author   = {Williams, Matthew O. and Kevrekidis, Ioannis G. and Rowley, Clarence W.},
  journal  = {J. Nonlinear Sci.},
  title    = {A data-driven approximation of the {K}oopman operator: extending dynamic mode decomposition},
  year     = {2015},
  issn     = {0938-8974,1432-1467},
  number   = {6},
  pages    = {1307--1346},
  volume   = {25},
  doi      = {10.1007/s00332-015-9258-5},
  mrnumber = {3415049},
}

@Article{Drmac2018,
  author   = {Drma\v{c}, Zlatko and Mezi\'{c}, Igor and Mohr, Ryan},
  journal  = {SIAM J. Sci. Comput.},
  title    = {Data driven modal decompositions: analysis and enhancements},
  year     = {2018},
  issn     = {1064-8275,1095-7197},
  number   = {4},
  pages    = {A2253--A2285},
  volume   = {40},
  doi      = {10.1137/17M1144155},
  fjournal = {SIAM Journal on Scientific Computing},
  mrnumber = {3828861},
}

@Article{Colbrook2023,
  author   = {Colbrook, Matthew J. and Ayton, Lorna J. and Sz\H{o}ke, M\'{a}t\'{e}},
  journal  = {J. Fluid Mech.},
  title    = {Residual dynamic mode decomposition: robust and verified {K}oopmanism},
  year     = {2023},
  issn     = {0022-1120,1469-7645},
  pages    = {Paper No. A21, 38},
  volume   = {955},
  doi      = {10.1017/jfm.2022.1052},
  fjournal = {Journal of Fluid Mechanics},
  mrnumber = {4534865},
}

@Article{Patel2015,
  author        = {Patel, Raajen and Goldstein, Thomas A. and Dyer, Eva L. and Mirhoseini, Azalia and Baraniuk, Richard G.},
  journal       = {ar{X}iv preprint},
  title         = {{oASIS}: Adaptive Column Sampling for Kernel Matrix Approximation},
  year          = {2015},
  archiveprefix = {arXiv},
  eprint        = {1505.05208},
  primaryclass  = {stat.ML},
  publisher     = {arXiv},
}

@Article{Litzinger2018,
  author   = {Litzinger, Florian and Boninsegna, Lorenzo and Wu, Hao and N{\"u}ske, Feliks and Patel, Raajen and Baraniuk, Richard and No{\'e}},
  journal  = {J. Chem. Theory Comput.},
  title    = {Rapid Calculation of Molecular Kinetics Using Compressed Sensing},
  year     = {2018},
  number   = {5},
  pages    = {2771--2783},
  volume   = {14},
  doi      = {10.1021/acs.jctc.8b00089},
  fjournal = {Journal of Chemical Theory and Computation},
}

@InCollection{Colbrook2024,
  author    = {Matthew J. Colbrook},
  booktitle = {Numerical Analysis Meets Machine Learning},
  publisher = {Elsevier},
  title     = {Chapter 4 - {T}he multiverse of dynamic mode decomposition algorithms},
  year      = {2024},
  isbn      = {9780443239847},
  pages     = {127--230},
  series    = {Handbook of Numerical Analysis},
  volume    = {25},
  doi       = {10.1016/bs.hna.2024.05.004},
  issn      = {1570-8659},
}

@InProceedings{Kawahara2016,
  author    = {Kawahara, Yoshinobu},
  booktitle = {Advances in Neural Information Processing Systems},
  title     = {Dynamic Mode Decomposition with Reproducing Kernels for {K}oopman Spectral Analysis},
  year      = {2016},
  publisher = {Curran Associates, Inc.},
  volume    = {29},
}

@Article{Askham2018,
  author   = {Askham, Travis and Kutz, J. Nathan},
  journal  = {SIAM J. Appl. Dyn. Syst.},
  title    = {Variable projection methods for an optimized dynamic mode decomposition},
  year     = {2018},
  issn     = {1536-0040},
  number   = {1},
  pages    = {380--416},
  volume   = {17},
  doi      = {10.1137/M1124176},
  fjournal = {SIAM Journal on Applied Dynamical Systems},
  mrnumber = {3796201},
}

@InProceedings{Schoelkopf2000,
  author    = {Sch\"{o}lkopf, Bernhard},
  booktitle = {Advances in Neural Information Processing Systems},
  title     = {The Kernel Trick for Distances},
  year      = {2000},
  publisher = {MIT Press},
  volume    = {13},
}

@InProceedings{Baudat2001,
  author    = {Baudat, G. and Anouar, F.},
  booktitle = {Proceedings of the International Joint Conference on Neural Networks, IEEE},
  title     = {Kernel-based methods and function approximation},
  year      = {2001},
  pages     = {1244--1249},
  volume    = {2},
  doi       = {10.1109/IJCNN.2001.939539},
}

@InProceedings{Williams2000,
  author    = {Williams, Christopher and Seeger, Matthias},
  booktitle = {Advances in Neural Information Processing Systems},
  title     = {Using the Nystr\"{o}m Method to Speed Up Kernel Machines},
  year      = {2000},
  publisher = {MIT Press},
  volume    = {13},
}

@Book{Zhang2005,
  editor    = {Zhang, Fuzhen},
  publisher = {Springer-Verlag, New York},
  title     = {The {S}chur complement and its applications},
  year      = {2005},
  isbn      = {0-387-24271-6},
  series    = {Numerical Methods and Algorithms},
  volume    = {4},
  doi       = {10.1007/b105056},
  mrnumber  = {2160825},
  pages     = {xvi+295},
}

@Article{Schmid2010,
  author   = {Schmid, Peter J.},
  journal  = {J. Fluid Mech.},
  title    = {Dynamic mode decomposition of numerical and experimental data},
  year     = {2010},
  issn     = {0022-1120,1469-7645},
  pages    = {5--28},
  volume   = {656},
  doi      = {10.1017/S0022112010001217},
  fjournal = {Journal of Fluid Mechanics},
  mrnumber = {2669948},
}

@Article{Tu2014,
  author   = {Tu, Jonathan H. and Rowley, Clarence W. and Luchtenburg, Dirk M. and Brunton, Steven L. and Kutz, J. Nathan},
  journal  = {J. Comput. Dyn.},
  title    = {On dynamic mode decomposition: theory and applications},
  year     = {2014},
  issn     = {2158-2491,2158-2505},
  number   = {2},
  pages    = {391--421},
  volume   = {1},
  doi      = {10.3934/jcd.2014.1.391},
  fjournal = {Journal of Computational Dynamics},
  mrnumber = {3415261},
}

@Book{Kutz2016,
  author    = {Kutz, J. Nathan and Brunton, Steven L. and Brunton, Bingni W. and Proctor, Joshua L.},
  publisher = {Society for Industrial and Applied Mathematics},
  title     = {Dynamic mode decomposition},
  year      = {2016},
  address   = {Philadelphia, PA},
  isbn      = {978-1-611974-49-2},
  doi       = {10.1137/1.9781611974508},
  mrnumber  = {3602007},
}

@Book{Higham2002,
  author    = {Higham, Nicholas J.},
  publisher = {Society for Industrial and Applied Mathematics},
  title     = {Accuracy and Stability of Numerical Algorithms},
  year      = {2002},
  address   = {Philadelphia, PA, USA},
  edition   = {Second},
  isbn      = {0-89871-521-0},
  doi       = {10.1137/1.9780898718027},
  pages     = {xxx+680},
}
\end{document}